\documentclass[12pt,a4paper]{article}
\usepackage{latexsym,amsmath,amssymb}
\usepackage[active]{srcltx}
\date{}
\numberwithin{equation}{section}

\DeclareMathOperator{\co}{co}

\DeclareMathOperator{\inter}{int} 

\DeclareMathOperator{\aff}{aff}

\DeclareMathOperator{\clo}{cl}

\DeclareMathOperator{\sco}{sco}

\DeclareMathOperator{\dist}{dist}

\newcommand{\R}{\mathbb R}

\newcommand{\sbs}{\subset}

\newtheorem{theo}{Theorem}
\newtheorem{lem}{Lemma}

\newtheorem{remark}{Remark}

\newcommand{\rang}{\rangle}
\newcommand{\lang}{\langle}

\begin{document}

\title{Sphere covering by minimal number of caps and short closed sets
\thanks{{\it 1991 A M S Subject Classification.} 52A45
{\it Key words and phrases.} Sphere covering by closed sets.}}
\author{A. B. N\'emeth}
\maketitle

\begin{abstract}
A subset of the sphere is said short if it is contained in an open
hemisphere. A short closed set which is geodesically convex is called a cap.
The following theorem holds:
1. The minimal number of short closed sets covering the $n$-sphere is $n+2$.
2. If $n+2$ short closed sets cover the $n$-sphere then
(i) their intersection is empty;
(ii) the intersection of any proper subfamily of them is non-empty.
In the case of caps (i) and (ii) are also sufficient for the family to
be a covering of the sphere.
\end{abstract}

\section{Introduction and the main result}

Denote by $\R^{n+1}$ the $n+1$-dimensional Euclidean space endowed
with a Cartesian reference system,
with the scalar product $\lang\cdot,\cdot\rang$ and with the topology it generates.
 
Denote by $S^n$ the $n$-dimensional unit sphere in $\R^{n+1}.$

A subset of the sphere $S^n$  is said \emph{short} if it is contained in an open
hemisphere. 

The subset $C\sbs S^n$ is called {\it geodesically convex}
if together with any two of its points it contains the arc
of minimal length of the principal circle on $S^n$ through these
points. $S^n$ itself is a geodesically convex set. 

A short closed set which is geodesically convex is called a \emph{cap}. 

We use the notation $\co A$ for the convex hull of $A$ and the notation 
$\sco A$ for the geodesical convex hull  of $A\subset S^n$ (the union of the
geodesical lines with endpoints in $A$). Further $\dist(\cdot,\cdot)$ will denote
the geodesical distance of points.
Besides the standard notion of simplex we also use 
the notion of the spherical simplex $\Delta$
placed in the north hemisphere $S^+$ of $S^n$ such that their vertices are on the equator of $S^n$.
In this case
$\|\Delta\|=S^+$.

Our main result is:

\begin{theo}
\begin{enumerate}

\item The minimal number of short closed sets covering $S^n$ is $n+2$.

\item If a family  $F_1,...F_{n+2}$  of short closed sets covers $S^n$, then:

(i) $\cap_{i=1}^{n+2} F_i = \emptyset$;

(ii) $\cap_{i\not=j}F_i\not= \emptyset,\;\forall \; j=1,...,n+2$;

(iii) if $a_j\in \cap_{i\not=j}F_i$, then the vectors $a_1,...,a_{n+2}$
are the vertices of an $n+1$-simplex containing $0$ in its interior.

\end{enumerate}

If the sets $F_i$ are caps, then (i) and (ii) are also sufficient for the family to be a cover of $S^n$.
\end{theo}

 Let $\Delta$ be an $n+1$-dimensional simplex with
vertices in $S^n$ containing the origin in its interior.
Then the radial projection from $0$ of the closed $n$-dimensional
faces of $\Delta$ into $S^n$ furnishes $n+2$ caps
covering $S^n$ and satisfying  (i) and (ii).

A first version for caps of the above theorem is the content of the unpublished note \cite{nemeth2006}.

\begin{remark}
We mention the formal relation in case of caps with those in the \emph{Nerve Theorem} (\cite{hatcher2002} Corollary 4G3).
If we consider \emph{"open caps"} in place of caps, then the conclusion (ii) can be deduced from
the mentioned theorem. Moreover, the conclusion holds for a "good" open cover of the sphere too, 
i. e., an open cover with contractible members and contractible  finite intersections.
In our theorem the covering with caps has the properties of a "good" covering in this theorem:
the members of the covering together with their nonempty intersections are contractible, but
their members are closed, circumstance which seems to be rather sophisticated to be surmounted.
(Thanks are due to Imre B\'ar\'any, who mentioned me this possible connection.)
\end{remark}

We shall use in the proofs the following (spherical)
variant of Sperner's lemma
(considered for simplices by Ky Fan \cite{FA}):
\begin{lem}\label{sperner}
If a collection of closed sets $F_1,...,F_{n+1}$ in $S^n$
covers the spherical simplex engendered by the points $a_1,...a_{n+1}\in S^n$ and 
$\sco \{a_1,...,a_{i-1},a_{i+1},...,a_{n+1}\}\subset F_i,\; i=1,...,n+1$ then
$\cap_{i=1}^{n+1} F_i\not= \emptyset.$
\end{lem}

Our first goal is to present the proof for caps. (We mention that using the methods in
\cite{KN} and \cite{KN1} the proof can be carried out in a purely geometric way in contrast with the proof in
\cite{nemeth2006}, where we refer to the Sperner lemma.)

Using the variant for caps of the theorem and the Sperner lemma we prove then the
variant for short closed sets.
  
Except the usage of Lemma \ref{sperner}, our methods are elementary: 
they use repeatedly the induction with respect to the dimension.


\section{The proof of the theorem for caps}

1.
Consider $n+1$ caps $C^1,...,C^{n+1}$ on $S^n$. $C^i$ being a cap, can be separated
strictly by a hyperplane
$$H_i= \{x\in \R^{n+1}:\,\lang a_i,x\rang +\alpha_i =0\}$$
from the origin. We can suppose without loss of generality, that
the normals $a_i$ are linearly independent, since
by slightly moving them we can achieve this, without affecting
the geometrical picture.
If the normals $a_i$ are considered oriented toward $0$,
this strict separation means that $\alpha_i >0,\;i=1,...,n.$
The vectors $a_i,\,i=1,...,n+1$ engender a reference
system in $\R^{n+1}$. Let $x$ be a nonzero element of
the positive orthant of this reference system. Then,
for $t\geq 0$, one has
$\lang a_i,tx \rang \geq 0,\;\forall \,i=1,...,n+1.$

Hence, for each $t\geq 0$,
$tx$ will be a solution of the system

$$\lang a_i,y\rang +\alpha_i>0,\;i=1,...,n+1.$$
and thus
$$(*)\qquad tx\in \cap_{i=1}^{n+1} H_i^+,\;\forall \, t\geq 0$$
with
$$H_i^+=\{y \in R^{n+1};\,\lang a_i,y\rang +\alpha_i > 0\}.$$

Now, if $C^1,...,C^{n+1}$ covers $S^n$, then
so does the union  $\cup_{i=1}^{n+1} H_i^-$ of halfspaces
$$H_i^-=\{y \in R^{n+1};\,\lang a_i,y\rang +\alpha_i \leq 0\}.$$

Since $H_i^+$ is the complementary set of $H_i^-$
and $S^n\sbs \cup_{i=1}^{n+1} H_i^-$, the set
$\cap_{i=1}^{n+1} H_i^+$ must be inside $S^n$ and
hence bounded. But (*) shows that $tx$ with
$x\not= 0$ is in this set for any $t\geq 0$.
The obtained contradiction shows that the
family $C^1,...,C^{n+1}$ cannot cover $S^n$.

\begin{remark}

The proof of this item is also consequence
of the Lusternik-Schnirelmann theorem
\cite{LS} which asserts that if $S^n$
is covered by the closed sets $F_1,...,F_k$ with
$F_i\cap (-F_i)=\emptyset,\,i=1,...,k,$ then $k\geq n+2.$

\end{remark}

2. Let $C^1,...,C^{n+2}$ be caps covering $S^n$. 

(i) Then they cannot have a common
point $x$, since this case $-x$ cannot be covered by any
$C^i$. (No cap
can contain diametrically opposite points of $S^n$.)

Hence, condition (i) must hold.

(ii) To prove that $\cap_{j\not= i} C^j \not= \emptyset,
\;\forall \,i=1,...,n+2$ we proceed by induction.

For $S^1$, the circle, $C^i$ is an arc (containing
its endpoints) of length $< \pi$,
$i=1,2,3$. The arcs $C^1, C^2, C^3$ cover $S^1$.
Hence, they cannot have common points, and the endpoint
of each arc must be contained in exactly one of
the other two arcs. Hence, $C^i$ meets $C^j$ for
every $j\not= i.$ If $c_i\in C^j\cap C^k,\;
j\not= i\not= k\not=j$, then $c_1, \,c_2,\,c_3$ are tree
pairwise different points on the circle,
hence they are in general position and $0$ is
an interior point of the  triangle they span.

Suppose the assertions (ii) and (iii) hold for
$n-1$ and let us prove them for $n$.

Take $C^{n+2}$ and let $H$ be a hyperplane through
$0$ which does not meet $C^{n+2}.$ Then, $H$ determines
the  closed hemispheres $S^-$ and $S^+$. Suppose
that $C^{n+2}$ is placed inside $S^-$ (in the interior
of $S^-$ with respect the topology of $S^n$).
Hence, $C^1,...,C^{n+1}$ must cover $S^+$ and
denoting by $S^{n-1}$ the $n-1$-dimensional sphere $S^n\cap H$, 
these sets cover $S^{n-1}.$
Now, $D^i= C^i\cap S^{n-1},\;i=1,...,n+1$ are
caps in $S^{n-1}$
which cover this sphere. Thus, the induction hypothesis
works for these sets.

Take the points $d_i\in \cap_{j\not= i}D^j$. Then,
$d_1,...,d_{n+1}$ will be in general position and $0$ is an
interior point of the simplex they span. By
their definition, it follows that $d_k\in D^j, \;
\forall k\not= j$ and hence $d_1,...,d_{j-1},d_{j+1},
...,d_{n+1}\in D^j,\; j=1,2,...n+1. $

Consider the closed hemisphere $S^+$ to be endowed
with a spherical simplex structure $\Delta$
whose vertices are the points $d_1,...,d_{n+1}$
.

Since $C^1,...,C^{n+1}$ cover $S^+$, and
$d_1,...,d_{j-1},d_{j+1},
...,d_{n+1}\in D^j\subset C^j\cap S^+,\; j=1,2,...n+1 $,
Lemma \ref{sperner}
 can be applied to the spherical simplex $\Delta$, yielding
$$\cap_{j=1}^{n+1} C^j \supset\cap_{j=1}^{n+1} (S^+\cap C^j) \not=
\emptyset.$$

This shows that each collection of $n+1$ sets $C^j$ have
nonempty intersection and proves (ii) for $n$.

(If we prefer a purely geometric proof of this item, we can refer to the
spherical analogue of the results in \cite{KN1}.)

From the geometric picture is
obvious that two caps meet if and only
if their convex hulls meet. Hence, from the conditions
(i) and (ii) for the caps $C^i$, it follows that
these conditions hold also for
$A^i=\co C^i,\;i=1,...,n+2.$

Take 
$$a_i\in \cap_{j\not= i}A^j,\;i=1,...,n+2.$$

Let us show that
for an arbitrary $k\in N$, 
$$a_k\not\in
 \aff \{a_1,...,a_{k-1},a_{k+1},...,a_{n+2}\}.$$

 Assume the contrary. Denote
 $$H=\aff \{a_1,...,a_{k-1},a_{k+1},...,a_{n+2}\}.$$
Thus, $\dim H\leq n.$ The points $a_i$ are all in the manifold $H$. Denote
$$B^i=H\cap A^i.$$
Since $a_i\in \cap_{j\not= i}A^j$ and $a_i\in H$, it follows that
$$a_i\in \cap_{j\not= i}A^j\cap H=\cap_{j\not= i} B^j,\;\forall\,i.$$ 
This means that the family of convex compact
sets $\{B^j:\,j=1,...,n+2\}$ in $H$ possesses the property that every $n+1$ of its
elements have nonempty intersection.
Then, by Helly's theorem, they have a
common point. But this would
be a point of $\cap_{i=1}^{n+2} A^i $ too, which contradicts (ii) for the sets $A^i$.

Hence,
every $n+2$ points $c_i\in \cap_{j\not= i} C^j \sbs \cap_{j\not= i}
A^j,\; i=1,...,n+2$ are in general position. 
Since $c_1,...,c_{i-1},c_{i+1},...,c_{n+2} \in C^i$, it follows
that the open halfspace determined by the hyperplane they engender 
containing $0$ contains also the point $c_i$. 

This proves (iii).

Suppose that the caps $C^1,...,C^{n+2}$ posses the properties in (i) and (ii). 
Then, the method in the above proof yields that the points
 $$c_i\in \cap_{j\not= i} C^j,\;i=1,...,n+2$$ engender  an $n+1$-simplex with $0$ 
in its interior and
 $$c_1,...,c_{i-1},c_{i+1},...,c_{n+2} \in C^i,\; i=1,...,n+2.$$
The radial projections of the $n$-faces of this simplex into $S^n$
obviously cover $S^n$. The union of these projections
are contained in $\cup_{i=1}^{n+2} C^i$. 

This completes the proof for cups.


\section{The proof of the theorem for short closed sets}

We carry out the proof by induction. 

Consider  $n=1$ and suppose $F_1,F_2,F_3$ are short closed sets  covering $S^1$.

If $a\in \cap_{i=1}^3 F_i$, then by the above hypothesis $-a\not \in \cup_{i=1}^3 F_i$, which is impossible.
Hence, 
$$\cap_{i=1}^3 F_i= \emptyset$$ must hold.

 Denote $C_3= \sco F_3$,
then $C=\clo S^1\setminus C_3$ is a connected arc of $S^1$ covered by $F_1,F_2$.
One must have $C\cap F_i\not=\emptyset , i=1,2$, since if for instance $C\cap F_2= \emptyset,$
then it would follow that the closed sets $F_1$ and $C_3$, both of geodesical diameter $< \pi$
cover $S^1$, which is impossible.    
Since $C$ is connected and $C\cap F_i,\; i=1,2$ are closed sets 
in $C$ covering this set, $F_1\cap F_2\supset (C\cap F_1)\cap (C\cap F_2) \not= \emptyset$ must hold.

The geodesically convex sets $C_i=\sco F_i,\;i=1,2,3$ cover $S^1$, hence applying
the theorem for caps to $ a_j  \in \cap_{i\not= j} F_i \subset \cap_{i\not= j} C_i,\; j=1,2,3$,
we conclude that these points are in general position and
the simplex engendered by them must contain $0$ as an interior point.

Suppose that the assertions hold for $n-1$ and prove them for $n$.

Suppose that 
$$S^n \subset \cup_{i=1}^{n+2} F_i,\; F_i\;\textrm{short, closed},\; i=1,...,n+2.$$

The assertion (i) is a consequence of the theorem for caps applied to
$C_i =\sco F_i,\; i=1,...,n+2$ (or a consequence of the Lusternik Schnirelmann  theorem).

Suppose that $F_{n+2}$ is contained in the interior (with respect to the topology
of $S^n$) of the south hemisphere $S^-$ and denote by $S^{n-1}$
the equator of $S^n$. 

 Now $S^{n-1} \subset \cup_{i=1}^{n+1} F'_i$ with $F'_i=(S^{n-1}\cap F_i
),\; i=1,...,n+1,$
and we can apply the induction hypothesis for $S^{n-1}$ and the closed sets $F'_i,\; i=1,...,n+1.$
Since $C'_i=\sco F'_i,\; i=1,...,n+1$ cover $S^{n-1}$, and they are caps,
the theorem for caps applies and hence the points
\begin{equation*} 
a_j \in \cap_{i=1,i\not=j}^{n+1} C'_i,\;j=1,...,n+1
\end{equation*}
are in general position.

The closed sets
\begin{equation*}
A_i = C'_i\cup(F_i\cap S^+)= C'_i\cup (F_i\cap \inter S^+),\;i=1,...,n+1
\end{equation*}
cover $S^+$, the north hemisphere considered as a spherical simplex $\Delta$ engendered by $a_1,...,a_{n+1}$ $(\|\Delta\|=S^+$).
(Here $\inter S^+$ is the interior of $S^+$
in the space $S^n$.)
Further, 
$$\sco \{a_1,...,a_{k-1},a_{k+1},...,a_{n+1}\} \subset A_k,\; k=1,...,n+1.$$
Hence, we can apply Lemma \ref{sperner} to conclude that there exists a point $a$ in
$\cap_{i=1}^{n+1} A_i \not= \emptyset.$  

Since 
$$C'_i\cap (F_j\cap \inter S^+)=\emptyset,\;\forall \;i,\;j,$$
it follows that
$$a \in \cap_{i=1}^{n+1} A_i= \cap_{i=1}^{n+1} C'_i \cup \cap_{i+1}^{n+1} F_i\cap \inter S^+=\cap_{i+1}^{n+1} F_i\cap \inter S^+,$$
because $\cap_{i=1}^{n+1} C'_i= \emptyset$ by the induction hypothesis and the theorem for caps.
Thus,
$$a\in \cap_{i=1}^{n+1} F_i,$$
and we have condition (ii) fulfilled for $n$.

The condition (iii) follows from the theorem for caps applied to
$$C^i= \sco F^i,\;i=1,...,n+2.$$

\begin{remark}
If $S^1$ is covered by the closed sets $F_1, F_2, F_3$ with the property
$F_i\cap (-F_i)=\emptyset,\; i=1,2,3 $, then
$F_i\cap F_j \not= \emptyset$ $\forall \;i, j.$

Indeed, assume that $F_1\cap F_2= \emptyset$. Then $\dist (F_1,F_2)=\varepsilon >0.$ If $a_i\in F_i$
are the points in $F_i, \; i=1,2$ with $\dist (a_1,a_2)= \varepsilon,$ then the closed arc $C\subset S^1$
with the endpoints $a_1,\; a_2$ must be contained in $F_3$, and hence $-C\cap F_3= \emptyset$, and
then $-C$ must be covered by $F_1\cup F_2$. Since $-a_1 \in -C$ cannot be in $F_1$, it must be in $F_2$,
and $-a_2\in F_1$. Thus, $F_1\cap -C\not= \emptyset$ and $F_2\cap -C \not= \emptyset,$ while
the last two sets cover $-C$. Since $-C$ is connected and the respective sets are closed, they
must have a common point, contradicting the hypothesis $F_1\cap F_2= \emptyset$.

This way, we obtain (ii) fulfilled for $n=1$ for this more general case.
We claim that the conditions also hold for $n$, that is, if the closed sets $F_1,...,F_{n+2}$ with
$F_i\cap (-F_i)=\emptyset,\; i=1,...,n+2 $ cover $S^n$, then condition (ii) holds. (Condition
(i) is a consequence of the definition of the sets $F_i$.)
\end{remark}

\end{document}